\documentstyle[12pt]{article} 
\newcommand{\Vdash}{\mathrel{{\vrule height 6.9pt depth -0.1pt}\!\vdash}} 
\newcommand{\restriction}{|} 

\newcommand{\razy}{\!\cdot\!}

\newcommand{\rest}{{\mathord{\restriction}}} 
 
\newcommand{\nor}{{\rm {\bf nor}}\/}

\newcommand{\lh}{{\rm lh}\/} 
\newcommand{\dom}{{\rm dom}} 
\newcommand{\rng}{{\rm rng}}
\renewcommand{\root}{{\rm root}\/} 
\newcommand{\suc}{{\rm succ}} 
\newcommand{\QED}{\hspace{0.2in}\vrule width 6pt height 6pt depth 0pt 
\vspace{0.1in}} 

\newcommand{\forces}{\Vdash} 
\newcommand{\R}{{\cal R}} 
\newcommand{\Q}{{\bf Q}} 
\newcommand{\p}{{\bf P}} 
\newcommand{\C}{{\bf C}} 
\newcommand{\D}{{\bf D}}  
\newcommand{\B}{{\bf B}}  
\renewcommand{\L}{{\bf L}} 
\newcommand{\can}{2^{\textstyle \omega}} 
\newcommand{\baire}{\omega^{\textstyle \omega}} 
 
\newcommand{\V}{{\bf V}} 

\newcommand{\lesdot}{\mathrel{\mathord{<}\!\!\raise 
0.8 pt\hbox{$\scriptstyle\circ$}}} 
\newcommand{\Proof}{{\sc Proof} \hspace{0.2in}} 

\newtheorem{theorem}{Theorem}[section] 
\newtheorem{lemma}[theorem]{Lemma} 
\newtheorem{proposition}[theorem]{Proposition} 
\newtheorem{corollary}[theorem]{Corollary} 
\newtheorem{problem}[theorem]{Problem}

\title{Examples for Souslin Forcing\thanks{The authors would like to thank Israel
Academy of Sciences BRF and MSRI for partial support.}} 
\author{{\bf Haim Judah}\\
Dept. of Mathematics and Computer Science\\
Bar-Ilan University\\
52900 Ramat-Gan, Israel
\and
{\bf Andrzej Ros\l anowski}\\
Dept. of Mathematics and Computer Science\\
Bar-Ilan University\\
52900 Ramat-Gan, Israel\\
and\\
Mathematical Institute of Wroclaw University\\
50384 Wroclaw, Poland
\and
{\bf Saharon Shelah}\\
Institute of Mathematics\\
The Hebrew University of Jerusalem\\
Jerusalem, Israel\\
and\\
Department of Mathematics\\
Rutgers University
Rutgers, USA}
 
\date\today 

\begin{document} 
\maketitle 

\section{Introduction} 

In this paper we continue with study of forcing notions having a simple 
definition. We began this study in [JS1] and [JS2]. In [BJ] we gave more 
results about Souslin forcing notions and in this paper we will give some
examples of Souslin forcing notions answering a question of [JS1] and a
question of H.Woodin. 

A forcing notion $\p$ is Souslin if $\p\subseteq\R$ is a $\Sigma_{1}^{1}$-set,
$\{(p,q): p\leq_{\p} q\}$ is a $\Sigma_{1}^{1}$-set and $\{(p,q): p \mbox{ is
incompatible with } q\}$ is a $\Sigma_{1}^{1}$-set. 

More information on Souslin forcing notion can be found in [JS1]. A related 
work is [BJ]. In [JS1] we prove that if $\p$ is Souslin ccc and $\Q$ is any
forcing notion then $\V^{\Q}\models$``$\p$ satisfies ccc''. A natural question
was: does ``$\p$ is Souslin ccc'' imply ``$\p$ has Knaster property''. 
Recall that $\p$ satisfies Knaster property if and only if 
\[(\forall A\!\in\![\p]^{\omega_{1}})(\exists B\!\in\![A]^{\omega_{1}})
(\forall p,q\in B)(p \mbox{ is compatible with } q).\]

In the second section we will give a model where there is a ccc Souslin 
forcing which does not satisfy the Knaster condition. Recall that under the
assumption of {\bf MA} every ccc notion of forcing has Knaster property.

Many simple forcing notions $\p$ satisfy the following condition:
\[\forces_{\p}\mbox{``}\hat{\p}\mbox{ is } \sigma\mbox{-centered''.}\] 
This property is connected with the homogeneity of the forcing notion.
The example of a totally nonhomogeneous Souslin forcing will be
constructed in the third section.

In the next section we present a model where there is a
$\sigma$-linked not $\sigma$-centered Souslin forcing such that all
its small subsets are $\sigma$-centered but Martin Axiom fails for
this order.

In section 5 we will give an example of a $\sigma$-centered Souslin forcing
notion and a model of the negation of {\bf CH} in which the union of less than
continuum meager subsets of $\R$ is meager but Martin Axiom fails for this 
notion of forcing.

In the last session of the MSRI Workshop on the continuum (October 1989)
H.Woodin asked if ``$\p$ has a simple definition and does not satisfy
ccc''implies that  
there exists a perfect set of mutually incompatible conditions. Clearly 
Mathias forcing satisfies such a requirement. In section 6 we will find a
Souslin forcing which is Proper but not ccc that does not contain a perfect
set of mutually incompatible conditions.

The last section will show that ccc $\Sigma^1_2$-notions of forcing
may not be indestructible ccc.

Our notation is standard and derived from [Je]. There is one exception,
however. We write $p\leq q$ to say that $q$ is a stronger condition then $p$.

\section{On the Knaster condition}

In this section we will build a Souslin forcing satisfying the countable chain
condition but which fails the Knaster condition.

Fix a sequence $<\sigma_{i}: i\in\omega>$ of functions from $\omega$ into 
$\omega$ such that  
\begin{description}
\item[{\bf (*)}] if $N<\omega$, $\phi_{i}:N\!\longrightarrow\!\omega$ (for 
$i<N$) then there are distinct $n_{0},\ldots,n_{N-1}$ such that
\[(\forall i,j_{0},j_{1}<N)(\phi_{i}(j_{0})=j_{1}\Rightarrow \sigma_{i}
(n_{j_{0}})=n_{j_{1}}).\]
\end{description}
Note that there exists a sequence $<\sigma_{i}: i\in\omega>$ satisfying 
{\bf (*)}:
\begin{quotation}
\noindent Suppose we have defined $\sigma_{i}\rest
m_{0}:m_{0}\!\longrightarrow\! m_{0}$ 
for $i<m_{0}$. We want to ensure {\bf (*)} for $n_{0}+1, \phi_{i}\ (i\leq 
n_{0})$. Define $\sigma_{i}(m_{0}+j_{0})=m_{0}+\phi_{i}(j_{0})$ for $i,j_{0}
\leq n_{0}$. Take large $m_{1}$ and extend all $\sigma_{i}$ ($i\leq n_{0}$) on
$m_{1}$ in such a way that $\rng(\sigma_{i})\subseteq m_{1}$. 
\end{quotation}
Next we define functions $f_{i}:\baire\longrightarrow\baire$ for $i\in\omega$ 
by
\[f_{i}(x)(k) = \left\{\begin{array}{ll}
                x(k)               & \mbox{if $k<i$}\\
                \sigma_{i}(x(k))   & \mbox{otherwise}
                \end{array}
                \right. \]      
Clearly all functions $f_{i}$ are continuous. Put $F(x)=\{f_{i}(x): i\in
\omega\}$ for $x\in\baire$. 

\begin{lemma}
\label{singletons}
\hspace{0.15in} Suppose that $x_{\alpha},y_{\alpha}\in\baire$ are such that
there is no repetition in $\{x_{\alpha},y_{\alpha}:\alpha\in\omega_{1}\}$.
Then there exists $A\in [\omega_{1}]^{\omega_{1}}$ such that
\[(\forall \alpha,\beta\in A)(\alpha<\beta \Rightarrow x_{\alpha}\not\in 
F(y_{\beta})).\]
\end{lemma}

\Proof For $\alpha<\omega_{1}$ let $n_{\alpha}=\min\{n:x_{\alpha}(n)\neq 
y_{\alpha}(n)\}$. We find $n\in\omega,\mbox{ and }s,t\in\omega^{\textstyle n}$
such that the set $A_{0}=\{\alpha<\omega_{1}: n=n_{\alpha}+1\ \&\ x_{\alpha}
\rest n=s \ \&\ y_{\alpha}\rest n=t\}$ is stationary in $\omega_{1}$. Clearly
$s\neq t$ and $s\rest (n-1)=t\rest (n-1)$. Thus $\alpha,\beta\in A_{0}$ and 
$x_{\alpha}\in F(y_{\beta})$ imply
$x_{\alpha}\in\{f_{i}(y_{\beta}):i\leq n\}$. 
Consequently the set $\{\alpha\in A_{0}\cap\beta: x_{\alpha}\in F(y_{\beta})
\}$ is finite for each $\beta\in A_{0}$.

We define the regresive function $\psi:A_{0}\longrightarrow\omega_{1}$ by 
$\psi(\beta)=\max\{\alpha\in A_{0}\cap\beta: x_{\alpha}\in F(y_{\beta})\}$ 
(with the convention that $\max\emptyset=0$). By Fodor's lemma there are 
$\gamma<\omega_{1}$ and a stationary set $A_{1}\subseteq A_{0}$ such that
$\psi(\beta)=\gamma$ for all $\beta\in A_{1}$. Put $A=A_{1}\backslash(\gamma
+1)$. Now, if $\alpha,\beta\in A$, $\alpha<\beta$ then $\psi(\beta)<\alpha$ 
and hence $x_{\alpha}\not\in F(y_{\beta})$. \QED

\begin{lemma}
\label{sets}
\hspace{0.15in} Suppose that $\{W_{\alpha}:\alpha<\omega_{1}\}$ is a family 
of disjoint finite subsets of $\baire$. Then there exist $\beta<\omega_{1}$
and an infinite set $A\subseteq\beta$ such that
\[(\forall \alpha\!\in\! A)(\forall x\!\in\! W_{\alpha})(\forall y\!\in\! 
W_{\beta})(x\not\in F(y)).\]
\end{lemma}

\Proof We may assume that all sets $W_{\alpha}$ are of the same cardinality,
say $\rest W\rest = n$ for $\alpha<\omega_{1}$. For $\alpha=\lambda+k$, where
$\lambda<\omega_{1}$ is a limit ordinal and $k\in\omega$ we define $X_{\alpha}
=W_{\lambda+2k}$ and $Y_{\alpha}=W_{\lambda+2k+1}$. Let $X_{\alpha}=\{
x_{i}^{\alpha}:i<n\}$, $Y_{\alpha}=\{y_{i}^{\alpha}:i<n\}$. Choose by the 
induction on $l=l_{1}\razy n+l_{2}<n^{2}$, $l_{1},l_{2}<n$ uncountable sets
$A_{l}\subseteq\omega_{1}$ satisfying
\begin{itemize}
\item $A_{l+1}\subseteq A_{l}$ and
\item if $l=l_{1}\razy n+l_{2}$, $l_{1},l_{2}<n$, $\alpha,\beta\in A_{l}$ and
$\alpha<\beta$ then $x^{\alpha}_{l_{1}}\not\in F(y^{\beta}_{l_{2}})$.
\end{itemize}
Since there is no repetition in $\{x^{\alpha}_{l_{1}},y^{\alpha}_{l_{2}}:
\alpha\in A_{l-1}\}$ we may apply lemma~\ref{singletons} to get $A_{l}$ from
$A_{l-1}$.

Consider $A_{n^{2}-1}$. Choose $\beta_{0}\in A_{n^{2}-1}$ such that the set 
\[A=\{\lambda+2k<\beta_{0}: \lambda+k\in A_{n^{2}-1}\ \& \ k\in\omega \ \& \ 
\lambda \mbox{ is a limit ordinal }\}\]
is infinite. Let $\beta_{0}=\lambda_{0}+k_{0}$ where $k_{0}\in\omega$ and
$\lambda_{0}$ is limit. Put $\beta=\lambda_{0}+2k_{0}+1$. Since $\beta_{0}<
\beta$ we have $A\subseteq\beta$. Suppose $\alpha=\lambda+2k\in A$. Let $x\in
W_{\alpha}$, $y\in W_{\beta}$. Then $\lambda+k\in A_{n^{2}_1}$, $W_{\alpha}=
X_{\lambda+k}$ and $W_{\beta}=Y_{\lambda_{0}+k_{0}}=Y_{\beta_{0}}$. Thus for
some $l_{1},l_{2}<n$ we have $x=x^{\lambda+k}_{l_{1}}$ and $y=y^{\beta_{0}}_{
l_{2}}$. Since $\lambda+k,\beta_{0}\in A_{l_{1}\cdot n+l_{2}}$, $\lambda+k<
\beta_{0}$ we get $x\not\in F(y)$. The lemma is proved. \QED

Let relations $R_{i}$ on $\omega^{\textstyle <\omega}$ be defined by
\begin{quotation}
$sR_{i}t$ if and only if

$i<\rest s\rest =\rest t \rest$, $s\rest i=t\rest i$ and $(\forall l\!\in\!
[i,\rest s\rest))(s(l)=\sigma_{i}(t(l)))$.
\end{quotation}
Note that if $x,y\in\baire$ are such that $(\forall n>i)(x\rest n R_{i} y\rest
n)$ then $x=f_{i}(y)$.

We define the following forcing notion $\Q$. A member $q$ of $\Q$ is a finite 
function such that:
\begin{description}
\item[$\alpha$) ] $\dom(q)\in [\omega_{1}]^{\textstyle <\omega}$, $\rng(q)
\subseteq\omega^{\textstyle <\omega}$,
\item[$\beta$) ] $(\forall \alpha,\beta\in\dom(q))(\alpha\neq\beta\Rightarrow
q(\alpha)\neq q(\beta))$,
\item[$\gamma$) ] there is $n(q)\in\omega$ such that $q(\alpha)\in\omega^{
\textstyle n(q)}$ for all $\alpha\in\dom(q)$.
\end{description}
The order is defined as follows:
\begin{quote}
$q\leq p$ if and only if
\end{quote}
\begin{enumerate}
\item $\dom(q)\subseteq\dom(p)$ and
\item $(\forall \alpha\in\dom(q))(q(\alpha)\subseteq p(\alpha))$ and
\item if $\alpha,\beta\in\dom(q)$, $\alpha<\beta$, $i<n(q)$ and $q(\alpha)
R_{i}q(\beta)$ then $p(\alpha)R_{i}p(\beta)$.
\end{enumerate}

\begin{lemma}
\label{ccc}
\hspace{0.15in} $\Q$ satisfies ccc.
\end{lemma}

\Proof Suppose $\{q_{\alpha}: \alpha<\omega_{1}\}\subseteq\Q$. We find $\gamma
<\omega_{1}$ and $A\in [\omega_{1}]^{\omega_{1}}$ such that for each $\alpha,
\beta\in A$, $\alpha<\beta$ we have
\begin{itemize}
\item $n(q_{\alpha})=n(q_{\beta})$,
\item $\dom(q_{\alpha})\cap\gamma=\dom(q_{\beta})\cap\gamma$, $(\dom(q_{
\alpha})\backslash\gamma)\cap(\dom(q_{\beta})\backslash\gamma)=\emptyset$, 
\item $q_{\alpha}\rest(\dom(q_{\alpha})\cap\gamma)=q_{\beta}\rest(\dom
(q_{\beta})\cap\gamma)$.
\end{itemize}
Suppose $\alpha,\beta\in A$. Clearly $\bar{q}=q_{\alpha}\cup q_{\beta}$ is a
function. The only problem is that there may exist $\gamma_{0}\in\dom(q_{
\alpha})$ and $\gamma_{1}\in\dom(q_{\beta})$ such that $q_{\alpha}(\gamma_{0}
)=q_{\beta}(\gamma_{1})$. Therefore to get a condition above both $q_{\alpha}$
and $q_{\beta}$ we have to extend all $\bar{q}(\gamma)$. Let $\dom(\bar{q})=
\{\gamma_{j}:j<N\}$ be an increasing enumeration. For $i<n(q_{\alpha})$ choose
$\phi_{i}:N\longrightarrow\omega$ such that
\begin{quotation}
\noindent if $j_{1}<j_{0}<N$ and
$\bar{q}(\gamma_{j_{1}})R_{i}\bar{q}(\gamma_{j_{0}})$ 
and either $\gamma_{j_{0}},\gamma_{j_{1}}\in\dom(q_{\alpha})$ or 
$\gamma_{j_{0}},\gamma_{j_{1}}\in\dom(q_{\beta})$

\noindent then $\phi_{i}(j_{0})=j_{1}$.
\end{quotation}
Note that $\bar{q}(\gamma_{j_{1}})R_{i}\bar{q}(\gamma_{j_{0}})$, $\bar{q}(
\gamma_{j_{2}})R_{i}\bar{q}(\gamma_{j_{0}})$ imply $\bar{q}(\gamma_{j_{1}})=
\bar{q}(\gamma_{j_{2}})$. Hence if $j_{2}<j_{1}<j_{0}$ are as above then
$\gamma_{j_{2}}\geq\gamma$ and consequently only one pair $(j_{1},j_{0})$ or
$(j_{2},j_{0})$ will be considered in the definition of $\phi_{i}$. Apply 
condition {\bf (*)} to find distinct $n_{0},\ldots,n_{N-1}$ such that
\[(\forall i<n)(q_{\alpha})(\forall j_{0},j_{1}<N)(\phi_{i}(j_{0})=j_{1}
\Rightarrow \sigma_{i}(n_{j_{0}})=n_{j_{1}}).\]
Put $q(\gamma_{j})=\bar{q}(\gamma_{j})\hat{\ }n_{j}$ for $j<N$. Clearly 
$q\in\Q$. Suppose $\gamma_{j_{1}},\gamma_{j_{0}}\in\dom(q_{\alpha})$,
$j_{1}<j_{0}$ and
$q_{\alpha}(\gamma_{j_{1}})R_{i}q_{\alpha}(\gamma_{j_{0}})$ for 
some $i<n(q_{\alpha})$. Then $\phi_{i}(j_{0})=j_{1}$ and hence
$\sigma_{i}(n_{j_{0}})=n_{j_{1}}$. Thus
$q(\gamma_{j_{1}})R_{i}q(\gamma_{j_{0}})$. It  
shows that $q_{\alpha}\leq q$. Similarly $q_{\beta}\leq q$. Thus we have 
proved that $\Q$ satisfies Knaster condition. \QED

Let $G\subseteq \Q$ be generic over $\V$. In $\V[G]$ we define $x_{\alpha}^{G}
=\bigcup\{q(\alpha):q\in G\ \&\ \alpha\in\dom(q)\}$ for $\alpha<\omega_{1}$.
Obviously each $x_{\alpha}^{G}$ is a sequence of integers. As in the proof of
lemma~\ref{ccc} we can show that for each $q\in\Q$ there is $p\geq q$ such 
that $n(p)=n(q)+1$. Consequently $x_{\alpha}^{G}\in\baire$ for every $\alpha<
\omega_{1}$. Moreover $x_{\alpha}^{G}\neq x_{\beta}^{G}$ for $\alpha<\beta<
\omega_{1}$ (recall that $q(\alpha)\neq q(\beta)$ for distinct $\alpha,\beta
\in\dom(q)$).

Note that if $\alpha<\beta$, $\alpha,\beta\in\dom(q)$, $q\in\Q$ and $i<n(q)$ 
then
\begin{quotation}
$q(\alpha)R_{i}q(\beta)$ implies $q\forces \dot{x}_{\alpha}=f_{i}(\dot{x}_{
\beta})$ and

$\neg q(\alpha)R_{i}q(\beta)$ implies $q\forces \dot{x}_{\alpha}\neq f_{i}(
\dot{x}_{\beta})$.
\end{quotation}

\begin{lemma}
\label{model}
\hspace{0.15in} Suppose $G\subseteq \Q$ is generic over $\V$. Then
\[\V[G]\models (\forall A\!\in\![\omega_{1}]^{\omega_{1}})(\exists \alpha\!,
\!\beta\!\in\! A)(\alpha<\beta \ \& \ x_{\alpha}^{G}\in F(x_{\beta}^{G})).\]
\end{lemma}

\Proof Let $\dot{A}$ be a $\Q$-name for an uncountable subset of $\omega_{1}$.
Given $p\in\Q$. Find $A_{0}\in [\omega_{1}]^{\omega_{1}}$ and $q_{\alpha}\geq
p$ for $\alpha\in A_{0}$ such that $\alpha\in\dom(q_{\alpha})$ and $q_{\alpha}
\forces\alpha\in\dot{A}$. We may assume that for each $\alpha,\beta\in A_{0}$
we have $n=n(q_{\alpha})=n(q_{\beta})$ and $q_{\alpha}(\alpha)=q_{\beta}
(\beta)$. Now we repeat the procedure of lemma~\ref{ccc} with one small 
change. We choose suitable $A_{1}\in [A_{0}]^{\omega_{1}}$, $\gamma<
\omega_{1}$ and we take $\alpha,\beta\in A_{1}$, $\alpha<\beta$. Defining 
integers $n_{0},\ldots,n_{N-1}$ we consider functions $\phi_{i}:N
\longrightarrow\omega$ (for $i<n$) as in \ref{ccc} and a function
$\phi_{n}:N 
\longrightarrow\omega$ such that $\phi_{n}(k)=l$, where $\alpha=\gamma_{l}$,
$\beta=\gamma_{k}$. Then we get a condition $q\in\Q$ above both $q_{\alpha}$
and $q_{\beta}$ and such that $\sigma_{n}(q(\beta)(n))=q(\alpha)(n)$.
Since $q(\beta)\rest n=q(\alpha)\rest n$ and $n(q)=n+1$ we have $q(\alpha)
R_{n}q(\beta)$ and consequently $q\forces \dot{x}_{\alpha}=f_{n}(\dot{x}_{
\beta})$. Since $q\geq q_{\alpha},q_{\beta}$ we get $q\forces\mbox{``}\alpha,
\beta\in\dot{A}\ \&\ \dot{x}_{\alpha}\in F(\dot{x}_{\beta})$''. \QED

Fix a Borel isomorphism 
$(\pi_{0},\pi_{1},\pi_{2}):\baire\longrightarrow (\baire)^{\omega}\times
2^{\textstyle \omega\times\omega}\times \baire$.
Thus if $x\in\baire$ then $\pi_{1}(x)$ is a relation on $\omega$ and $\pi_{0}
(x)$ is a sequence of reals. Let $\Gamma$ consists of all reals $x\in\baire$
such that
\begin{enumerate}
\item $(\forall n\neq m)(\pi_{0}(n)\neq\pi_{0}(m))$
\item $\pi_{1}(x)$ is a linear order on $\omega$
\item $\pi_{2}(x)\in A_{x}=\{\pi_{0}(n):n\in\omega\}$ and it is the 
$\pi_{1}(x)$-last element of $A_{x}$.
\end{enumerate}
Note that in 3 we think of $\pi_{1}(x)$ as an order on $A_{x}$. We define 
relations $<_{\Gamma}$ and $\equiv_{\Gamma}$ on $\Gamma$ by
\begin{quotation}
\noindent $x<_{\Gamma}y$ if and only if

\noindent $A_{x}$ is a proper $\pi_{1}(y)$-initial segment of $A_{y}$
and $\pi_{1}(y) \rest A_{x}=\pi_{1}(x)$.
\end{quotation}
\begin{quotation}
\noindent $x\equiv_{\Gamma}y$ if and only if 

\noindent $A_{x}=A_{y}$ and $\pi_{1}(y)=\pi_{1}(x)$ (we treat
$\pi_{1}(x),\pi_{1}(y)$ as  
orders on $A_{x},A_{y}$, respectively).
\end{quotation}
Clearly $\Gamma$ is a Borel subset of $\baire$, $<_{\Gamma}$ is a Borel 
transitive relation on $\Gamma$ and $\equiv_{\Gamma}$ is a Borel equivalence 
relation on $\Gamma$.

Now we define a forcing notion $\p_1$. Conditions in $\p_1$ are finite
subsets $p$ of $\Gamma$ such that
\begin{quote}
if $x,y\in p$, $x<_{\Gamma}y$ then $\pi_{2}(x)\not\in F(\pi_{2}(y))$.
\end{quote}
$\p_1$ is ordered by the inclusion.

\begin{lemma}
\label{P}
\hspace{0.15in} $\p_1$ is a ccc Souslin forcing.
\end{lemma}

\Proof $\p_1$ is Souslin since it can be easily coded as a Borel subset of 
$\baire$ in such a way that the order is Borel too. We have to show
that $\p_1$ 
satisfies the countable chain condition. First let us note some properties
of the incompability in $\p_1$. Suppose $p,q\in\p_1$ are incompatible. Clearly 
$p\backslash q$ and $q\backslash p$ are incompatible. If $x\in p$ and
$x\equiv_{\Gamma}x'$ then $(p\backslash\{x\})\cup\{x'\}$ and $q$ are 
incompatible. 

Suppose now that $\{p_{\alpha}:\alpha<\omega_{1}\}$ is an antichain in $\p_1$.
By the $\Delta$-lemma and by the above remarks we may assume that
\begin{description}
\item[1) ] $p_{\alpha}\cap p_{\beta}=\emptyset$ for $\alpha<\beta<\omega_{1}$
\item[2) ] if $x,x'\in \bigcup_{\alpha<\omega_{1}} p_{\alpha}$, $x\neq x'$ 
then $x\not\equiv_{\Gamma}x'$.
\end{description}
Note that if $p\in\p_1$ then the set 
$\{[y]_{\equiv_{\Gamma}}:y\in\Gamma \ \& \ (\exists x\in
p)(y<_{\Gamma}x)\}$ is countable. Hence, due to 2, we may assume that
\begin{description}
\item[3) ] $(\forall \alpha\!<\!\beta\!<\!\omega_{1})(\forall x\!\in\! 
p_{\alpha})(\forall y\!\in\! p_{\beta})(\neg y<_{\Gamma}x)$
\end{description}

{\sc Claim}: Let $d\in [\baire]^{\textstyle <\omega}$. Then $d=\{\pi_{2}(x):
x\in p_{\alpha}\}$ for at most countably many $\alpha<\omega_{1}$.

Indeed, assume not. Then we find $\beta<\omega_{1}$ such that $\{\pi_{2}(x):
x\in p_{\beta}\}=d$ and the set $B=\{\alpha<\beta: \{\pi_{2}(x):x\in 
p_{\alpha}\}=d\}$ is infinite. Note that if $x',x''<_{\Gamma}x$ and $\pi_{2}
(x')=\pi_{2}(x'')$ then $x'\equiv_{\Gamma}x''$. Hence if $x\in p_{\beta}$ then
for at most $\rest d\rest$ elements $x'$ of $\bigcup_{\alpha\in A} p_{\alpha}$
we have $x'<_{\Gamma}x$. Thus we find $\alpha\in A$ such that $(\forall
x'\!\in\! p_{\alpha})(\forall x\!\in\!p_{\beta})(\neg x'<_{\Gamma}x)$.
It follows from 3 that 
\[(\forall x\!\in\! p_{\beta})(\forall x'\!\in\! p_{\alpha})(\neg x
<_{\Gamma}x')\]
and hence conditions $p_{\alpha}$ and $p_{\beta}$ are compatible -- a 
contradiction.

Let $d_{\alpha}=\{\pi_{2}(x):x\in p_{\alpha}\}$. By the above claim we may 
assume that $d_{\alpha}\neq d_{\beta}$ for all $\alpha<\beta<\omega_{1}$.
Applying the $\Delta$-lemma we may assume that
\begin{description}
\item[4) ] $\{d_{\alpha}:\alpha<\omega_{1}\}$ forms a $\Delta$-system with
the root $d$.
\end{description}
Since the set $\bigcup_{w\in d} F(w)$ is countable w.l.o.g.
\begin{description}
\item[5) ] $(\forall \alpha\!<\!\omega_{1})(\forall v\in d_{\alpha}\backslash
d)(v\not\in\bigcup_{w\in d}F(w))$.
\end{description}
Apply lemma~\ref{sets} for the family $\{d_{\alpha}\backslash d:\alpha<
\omega_{1}\}$ to get $\beta<\omega_{1}$ and an infinite set $A\subseteq\beta$
such that
\begin{description}
\item[6) ] $(\forall \alpha\!\in\! A)(\forall v\in d_{\alpha}\backslash d)(
\forall w\in d_{\beta}\backslash d)(v\not\in F(w))$.
\end{description}
Let $y\in p_{\beta}$. As in the claim the set
\[\{x\in\bigcup_{\alpha\in A}p_{\alpha}: \pi_{2}(x)\in d \ \& \
x<_{\Gamma}y\}\] 
is finite. Consequently we find $\alpha\in A$ such that 
\begin{description}
\item[7) ] $(\forall x\!\in\! p_{\alpha})(\forall y\!\in\! p_{\beta})(\pi_{2}
(x)\in d \Rightarrow \neg x<_{\Gamma} y)$.
\end{description}
We claim that $p_{\alpha}$ and $p_{\beta}$ are compatible. Let $x\in 
p_{\alpha}$ and $y\in p_{\beta}$. By 7 we have that if $\pi_{2}(x)\in d$ then
$\neg x<_{\Gamma}y$. If $\pi_{2}(x)\not\in d$ and $\pi_{2}(y)\not\in
d$ then 6 applies and we get $\pi_{2}(x)\not\in F(\pi_{2}(y))$.
Finally if $\pi_{2}(x) 
\not\in d$ and $\pi_{2}(y)\in d$ then we use 5 to conclude that $\pi_{2}(x)
\not\in F(\pi_{2}(y))$. Hence $x\in p_{\alpha}$, $y\in p_{\beta}$ and 
$x<_{\Gamma}y$ imply $\pi_{2}(x)\not\in F(\pi_{2}(y))$. Consequently 
$p_{\alpha}\cup p_{\beta}\in\p_1$. \QED

\begin{lemma}
\label{notKnaster}
\hspace{0.15in} Assume that there exists a sequence $\{x_{\alpha}:\alpha<
\omega_{1}\}$ of elements of $\baire$ such that
\[(\forall A\!\in\![\omega_{1}]^{\omega_{1}})(\exists \alpha\!,\!\beta\!\in\! 
A)(\alpha<\beta \ \& \ x_{\alpha}\in F(x_{\beta})).\]
Then the forcing notion $\p_1$ does not satisfy the Knaster condition.
\end{lemma}

\Proof For $\alpha<\omega_{1}$ choose $y_{\alpha}\in\Gamma$ such that
\begin{itemize}
\item $A_{y_{\alpha}}=\{\pi_{0}(y_{\alpha})(n):n\in\omega\}=\{x_{\gamma}:\gamma
\leq\alpha\}$
\item $\pi_{1}(y_{\alpha})$ is the natural order on $A_{y_{\alpha}}$, 
$x_{\gamma}<_{\pi_{1}(y_{\alpha})}x_{\beta}$ iff $\gamma<\beta$.
\item $\pi_{2}(y_{\alpha})=x_{\alpha}$.
\end{itemize}
Let $p_{\alpha}=\{y_{\alpha}\}$ for $\alpha<\omega_{1}$. Then $\{p_{\alpha}:
\alpha<\omega_{1}\}$ does not have an uncountable subset of pairwise compatible
elements. \QED

Putting together lemmas~\ref{P}, \ref{notKnaster} and \ref{model} we get

\begin{theorem}
\label{thm}
\hspace{0.15in} It is consistent that there exists a ccc Souslin forcing notion
which does not satisfy the Knaster condition. \QED
\end{theorem}

It is not difficult to see that this example does not satisfy the following 
requirement:
\begin{quote}
``The generic object is encoded by a real''
\end{quote}
The next theorem says that also we can require such a condition. This answers
a question of J.Bagaria.

\begin{theorem}
\hspace{0.15in} It is consistent that there exists a ccc Souslin forcing 
notion $\Q$ such that $\forces_{\Q}\V[G]=\V[\dot{r}]$ for some $\Q$-name 
$\dot{r}$ for a real and $\Q$ does not satisfy the Knaster condition.
\end{theorem}

\Proof We follow the notation of the previous results. We work in the model 
of \ref{thm}. Let $\Q=\{(p,w):p\in\p \ \& \ w\in [\omega^{\textstyle 
<\omega}]^{\textstyle <\omega}\}$ be ordered by
\begin{quotation}
$(p,w)\leq (q,v)$ if and only if

$p\leq q$, $w\subseteq v$ and $x\rest n\not\in v\backslash w$ for every 
$x\in p$ and $n\in\omega$.
\end{quotation}
$\Q$ may be easily represented as a Souslin forcing notion (remember that
``for each $x$ in $p$'' is a quantification on natural numbers). Note that if
$p\leq q$, $p,q\in \p$ and $w\in [\omega^{\textstyle <\omega}]^{\textstyle
<\omega}$ then $(p,w)\leq (q,w)$. Hence $\Q$ satisfies the countable chain 
condition. If $\Q$ satisfied the Knaster condition then $\p$ would have 
satisfied it. We show that the $\Q$-generic object is encoded by a real. Let
$\dot{r}$ be a $\Q$-name for a subset of $\omega^{\textstyle <\omega}$ (a 
real) such that for any $\Q$-generic $G$ we have $\dot{r}^{G}=\bigcup\{w:
(\exists p)((p,w)\in G)\}$. Now in $\V[\dot{r}^{G}]$ define
\[H=\{(p,w)\in\Q: w\subseteq\dot{r}^{G}\ \&\ (\forall x\!\in\! p)(\forall n\!
\in\!\omega)(x\rest n\in\dot{r}^{G}\ \Leftrightarrow x\rest n\in w)\}.\]
Note that $H$ includes $G$ since $x\in p$, $x\rest n\not\in w$ imply $(p,w)
\forces x\rest n\not\in\dot{r}$. $H$ is a filter - suppose $(p_{0},w_{0}),
(p_{1},w_{1})\in H$. For each $x\in p_{0}\cup p_{1}$ we find $(p_{x},w_{x})
\in G$ such that $(p_{x},w_{x})\forces (\forall n\geq N)(x\rest n\not\in
\dot{r})$ for some $N$. If $x\not\in p_{x}$ we could take large $n$ and add
$x\rest n$ to $w_{x}$. Then we would have $(p_{x},w_{x})\leq (p_{x},w_{x}\cup
\{x\rest n\})$ and $(p_{x},w_{x}\cup\{x\rest n\})\forces x\rest n\in \dot{r}$. 
Thus $x\in p_{x}$ for all $x\in p_{0}\cup p_{1}$. Let $p=\bigcup_{x\in p_{0}
\cup p_{1}}p_{x}$, $w=\bigcup_{x\in p_{0}\cup p_{1}}w_{x}$. Then $(p,w)\in G
\subseteq H$ and $(p_{0},w_{0}),(p_{1},w_{1})\leq (p,w)$. Consequently $H=G$
and the theorem is proved. \QED

In the same time when the forcing notion $\p_1$ was constructed S.Todorcevic
found another example of this kind. 

Let $\cal F$ be the family of all converging sequences $s$ of real numbers 
such that $\lim s\notin s$. Todorcevic's forcing notion $\p_1^*$ consists
of finite subsets $p$ of $\cal F$ with property that
\[(\forall s,t\in p)(s\neq t\Rightarrow \lim s\notin t).\]
Todorcevic proved that $\p_1^*$ satisfies ccc and that if ${\bf b}=\omega_1$
then $\p_1^*$ does not have Knaster property (see [To]).

\section{A nonhomogeneous example}
\label{notsigma}
In this section we give an example of ccc Souslin forcing notion which
is very nonhomogeneous. Our forcing $\p_2$ will satisfy the following
property:
\begin{quotation}
{\em there is no $p\in \p_2$ such that}
$p\forces\mbox{`` $\hat{\p}_2\rest p$ is $\sigma$-centered''.}$
\end{quotation}
Recall that if $\Q$ is the Amoeba Algebra for Measure or the Measure
Algebra then $\forces_{\Q}\mbox{``$\hat{\Q}$ is $\sigma$-centered''}$
(see [BJ]). The Todorcevic example $\p_1^*$ has this property too. 
\begin{proposition}
\hspace{0.15in}$\forces_{\p_1^*}\mbox{``$\hat{\p_1^*}$ is
$\sigma$-centered''}$ 
\end{proposition}

\Proof For a rational number $d\in Q$ let $\phi_d:\p_1^*\longrightarrow
\p_1^*$ be the translation by $d$. Thus $\phi_d(p)=\{s+d: s\in p\}$. 
Note that $\phi_d$ is an automorphism of $\p_1^*$. Moreover if $p_1,p_2
\in\p_1^*$ then $\bigcup\{s-r:s\in p_1, r\in p_2\}$ is a nowhere dense set. 
Hence we find a rational $d$ such that
\[-d\notin\{a-b: a\in s\!\cup\!\{\lim s\}, b\in r\!\cup\!\{\lim r\},
s\in p_1, r\in p_2\}.\]
Then the conditions $\phi_d(p_1)$ and $p_2$ are compatible. Thus we have 
proved that for each $p\in\p_1^*$ the set $\{\phi_d (p): d\in Q\}$ is 
predense in $\p_1^*$. This implies that
\[\forces_{\p_1^*}\mbox{``}\hat{\p}_1^*=\bigcup_{d\in Q} \phi_d 
[\dot{\Gamma}]\mbox{ and each } \phi_d[\dot{\Gamma}]\mbox{ is 
centered'',}\]
where $\dot{\Gamma}$ is the canonical name for a generic filter. \QED

We do not know if
\[\forces_{\p_1}\mbox{``}\hat{\p}_1\mbox{ is }\sigma\mbox{-centered''.}\]
One can easily construct a ccc Souslin forcing $\p$ which does not force
that $\hat{\p}$ is $\sigma$-centered. An example of such a forcing
notion is the disjoint union of Cohen forcing and the measure algebra,
$\p=(\{0\}\times\C)\cup(\{1\}\times\B)$. In this order we have
$(0,\emptyset)\forces\mbox{``}\{1\}\times\hat{\B}\mbox{ is not
}\sigma$-centered''. But in this example we can find a dense set of
conditions $p\in\p$ such that
\[p\forces_{\p}\mbox{``}\hat{\p}\rest p=\{q\in\hat{\p}:q\geq p\}\mbox{
is }\sigma\mbox{-centered''}.\]

Define $T^*\subseteq\omega^{\textstyle <\omega}$,
$f,g:T^*\longrightarrow \omega$ in such a way that:
\begin{description}
\item[$(\alpha)$ ] $T^*$ is a tree,
\item[$(\beta)$ ] if $\eta\in T^*$ then $\suc_{T^*}(\eta)=f(\eta)$,
\item[$(\gamma)$ ] if $\lh\eta<\lh\nu$ or $\lh\nu=\lh\eta$ but
$(\exists k\!<\!\lh\eta)(\eta\rest k=\nu\rest k\ \&\ \eta(k)<\nu(k))$
then $f(\eta)<f(\nu)$,
\item[$(\delta)$ ] $g(\eta)>\rest T^*\cap \omega^{\textstyle
\lh\eta}\rest\cdot\prod\{f(\nu):f(\nu)<f(\eta)\}\cdot(100+\lh\eta)$ ,
\item[$(\epsilon)$ ] $f(\eta)> g(\eta)\cdot\prod\{2^{f(\nu)}:f(\nu)<
f(\eta)\}$
\end{description}

For $\eta\in T^*$ and and a set $A\subseteq\suc_{T^*}(\eta)=f(\eta)$
we define a norm of $A$:
\[\nor_\eta(A)=\frac{g(\eta)}{\rest f(\eta)\backslash A\rest}.\]
\begin{lemma}
\label{norm}
\hspace{0.15in}Suppose that $\eta\in T^*$ and $A_l\subseteq f(\eta)$ for
$l<m$. Let $\zeta=\min\{\nor_\eta(A_l):l<m\}$. Then

\noindent 1)\ \ \ $\nor_\eta(\bigcap_{l<m}A_l)\geq \zeta/m$,

\noindent 2)\ \ \ if $\zeta\geq 1$ and
$m\leq\prod\{2^{f(\nu)}:f(\nu)<f(\eta)\}$ then
$\bigcap_{l<m}A_l\neq\emptyset$. 
\end{lemma}

\Proof 1) Note that 
\[\rest f(\eta)\backslash\bigcap_{l<m}A_l\rest = \rest
\bigcup_{l<m}f(\eta)\backslash A_l\rest \leq \sum_{l<m}\rest
f(\eta)\backslash A_l\rest \leq m\cdot g(\eta)/\zeta\]
Hence 
\[\nor_{\eta}(\bigcap_{l<m}A_l) = \frac{g(\eta)}{\rest
f(\eta)\backslash\bigcap_{l<m}A_l\rest} \geq \zeta/m.\]

\noindent 2) Applying 1) we get $\nor_{\eta}(\bigcap_{l<m}A_l)\geq
1/m$. Hence 
\[\rest f(\eta)\backslash\bigcap_{l<m}A_l\rest\leq g(\eta)\cdot m \leq
g(\eta)\cdot\prod\{2^{f(\nu)}:f(\nu)<f(\eta)\}<f(\eta)\]
(the last inequality is guaranteed by condition ($\epsilon$)).
Consequently the set $\bigcap_{l<m}A_l$ is nonempty. \QED

Let $\p_2$ consists of all trees $T\subseteq T^*$ such that 
\[\lim_{n\rightarrow\infty}\min\{\nor_{\eta}(\suc_T(\eta)):\eta\in
T\cap \omega^{\textstyle n}\}=\infty.\]
The order is the inclusion.

Recall that a forcing notion $\Q$ is $\sigma$-$k$-linked if there
exist sets $R_n\subseteq\Q$ (for $n\in \omega$) such that
$\bigcup_{n\in \omega}R_n=\Q$ and each $R_n$ is $k$-linked (i.e. any
$k$ members of $R_n$ has a common upper boud in $\Q$).

\begin{proposition}
\hspace{0.15in}For every $k<\omega$ the forcing notion $\p_2$ is
$\sigma$-$k$-linked. 
\end{proposition}

\Proof Let $n\in \omega$ be such that for each $\eta\in T^*\cap
\omega^{\textstyle n}$
\[k<\prod\{2^f(\nu):f(\nu)<f(\eta)\}.\]
Note that the set
\[\{T\in\p_2:\lh(\root T)\geq n\ \&\ (\forall\eta\!\in\! T)(\root
T\subseteq \eta \Rightarrow \nor_{\eta}(\suc_T(\eta))\geq 1)\}\]
is dense in $\p_2$. For $\eta\in T^*$, $\lh\eta\geq n$ define
\[{\cal D}_\eta=\{T\!\in\!\p_2:\root T=\eta\ \&\ (\forall\\nu\!\in\!
T)(\eta\subseteq\nu \Rightarrow \nor_{\nu}(\suc_T(\nu))\geq
1)\}.\] 
Since $\bigcup\{{\cal D}_\eta: \lh\eta\geq n\}$ is dense in $\p_2$ it
is enough to show that each ${\cal D}_\eta$ is $k$-linked. Suppose
$T_0,\ldots,T_{k-1}\in{\cal D}_\eta$. Since
$k<\prod\{2^f(\nu):f(\nu)<f(\eta)\}$ we may apply lemma~\ref{norm}~2) to
conclude that if $\nu\in T=T_0\cap\ldots\cap T_{k-1}$,
$\nu\supseteq\eta$ then $\suc_T(\nu)\neq\emptyset$. By \ref{norm}~1)
we get $T\in\p_2$. \QED

For $\eta\in T^*$ we define the forcing notion $\Q_\eta$:
\[\begin{array}{ll}
\Q_\eta=\{t\subseteq T^*: & t\mbox{ is a finite tree of the
height } n\in \omega,\\
 &\root t =\eta \mbox{ and} \\
 &(\forall\nu\!\in\! t\!\cap\!\omega^{\textstyle
<n})(\eta\subseteq\nu\Rightarrow
\nor_{\nu}(\suc_t(\nu))\geq\lh\nu)\}.
\end{array}\] 
Since $\Q_\eta$ is countable and atomless it is isomorphic to Cohen
forcing $\C$. Let $\p=\prod\{\Q_{i,\eta}:i<\omega_1, \eta\in T^*\}$ be
the finite support product such that each $\Q_{i,\eta}$ is a copy of
$\Q_\eta$. 
\begin{theorem}
Let $G\subseteq \p$ be a generic filter over $\V$. Then, in $\V[G]$,
there is no $S\in\p_2$ such that
\[S\forces_{\p_2}\mbox{`` }\hat{\p}_2\rest S\mbox{ is
}\sigma\mbox{-centered''}.\] 
\end{theorem}

\Proof We work in $\V[G]$.

\noindent Assume $S\forces$``$\hat{\p}_2\rest S$ is
$\sigma$-centered''. Let $\dot{\R}_n$ ($n\in \omega$) be $\p_2$-names
for subsets of $\p_2$ such that 
\[S\forces_{\p_2}\mbox{`` }\hat{\p}_2\rest S\subseteq\bigcup_{n\in
\omega}\dot{R}_n \ \&\ \mbox{each } \dot{R}_n \mbox{ is directed''}.\]
Take $n\in \omega$ such that 
\[(\forall\eta\!\in\! S)(n\leq\lh\eta\Rightarrow
\nor_{\eta}(\suc_S(\eta)) \geq 1).\]
Fix any $\eta\in S\cap \omega^{\textstyle n}$ and choose
$l,m\in\suc_S(\eta)$, $l<m$. For $i<\omega_1$ put 
\[T_i=\bigcup\{t:\{((i,\eta\hat{\ }m),t)\}\in G\}.\]
Each $T_i$ is the tree added by $G\cap\Q_{i,\eta\hat{\ }m}$ and it is
an element of $\p_2$. Moreover $\root T_i=\eta\hat{\ }m$ and for each
$\nu\in T_i$ if $\eta\hat{\ }m\subseteq\nu$ then
$\nor_\nu(\suc_{T_i}(\nu))\geq \lh\nu$. Hence, by lemma~\ref{norm},
$T_i\cap S\in\p_2$ for each $i\in \omega_1$. 

\noindent Now we work in $\V$.

\noindent We find $p_i,\dot{S}_i,\eta_i,n_i$ such that for each $i\in
\omega_1$: 
\begin{itemize}
\item $p_i\in\p, n_i\in \omega,\eta_i\in T^*$ and $\dot{S}_i$ is a
$\p$-name for a member of $\p_2$,
\item $\forces_\p (\forall\nu\!\in\!\dot{S}_i)(\root\dot{S}_i
\subseteq\nu \Rightarrow \nor_\nu(\suc_{\dot{S}_i}(\nu))\geq 1)$,
\item $p_i\forces_\p\mbox{``}\eta\hat{\
}l\subseteq\eta_i=\root\dot{S}_i\ \&\ 
\dot{S}_i\forces_{\p_2}\dot{T}_i\in\dot{R}_{n_i}\mbox{''}$, 
\item $(i,\eta\hat{\ }m)\in\dom p_i$.
\end{itemize}
Next we find a set $I\in[\omega_1]^{\textstyle \omega_1}$ such that
$\{\dom p_i:i\in I\}$ forms a $\Delta$-system with the root $d$ and
for each $i\in I$:
\begin{itemize}
\item $\eta_i=\eta^*$ and $n_i=n^*$,
\item $p_i\rest d=p^*$,
\item $p_i(i,\eta\hat{\ }m)=t$, $(i,\eta\hat{\ }m)\notin d$.
\end{itemize}
Let $n^\#$ be the height of the tree $t$. Clearly we may assume that
$n^\#>\lh\eta^*$. Fix an enumeration $\{\rho_k:k<k^\#\}$ of $t\cap
\omega^{\textstyle n^\#}$. Put 
\[H=\{(a_k:k<k^\#): a_k\subseteq f(\rho_k)\ \&\ \nor_{\rho_k}(a_k)\geq
n^\#\}.\] 
Choose distinct $i_{\bar{\bf a}}\in I$ for $\bar{\bf a}\in H$. We define a
condition $q\in\p$ extending all $p_{i_{\bar{\bf a}}}$ ($\bar{\bf a}\in H$):
\begin{quotation}
\noindent $\dom q = \bigcup\{\dom p_{i_{\bar{\bf a}}}: \bar{\bf a}\in H\}$;

\noindent if $(i,\nu)\in\dom p_{i_{\bar{\bf a}}}$,
$(i,\nu)\neq(i_{\bar{\bf a}},\eta\hat{\ }m)$ then
$q(i,\nu)=p_{i_{\bar{\bf a}}}(i,\nu)$;

\noindent $q(i_{\bar{\bf a}},\eta\hat{\ }m)=t\cup\{\rho_k\hat{\ }c: k<k^\#,
c\in\bar{\bf a}(k)\}$. 
\end{quotation}
Now we take $r\geq q$ such that $r$ decides all
$\dot{S}_{i_{\bar{\bf a}}}\rest (n^\#+1)$. Thus we have finite trees
$s_{\bar{\bf a}}$ (for $\bar{\bf a}\in H$) such that $r\forces_\p
\dot{S}_{i_{\bar{\bf a}}}\rest (n^\#+1)=s_{\bar{\bf a}}$.
\smallskip

\noindent {\sc Claim:}\hspace{0.15in}: There exists $H'\subseteq H$
such that

\noindent {\em (i)}\hspace{0.15in} $\bigcap_{\bar{\bf a}\in
H'}(s_{\bar{\bf a}}\cap \omega^{\textstyle n^\#+1})\neq\emptyset$ and

\noindent {\em (ii)}\hspace{0.15in} for each $k<k^\#$ the set
$\bigcap\{\bar{\bf a}(k): \bar{\bf a}\in H'\}$ is empty.
\smallskip

Indeed, let $H_\rho=\{\bar{\bf a}\in H: \rho\in s_{\bar{\bf a}}\}$ for
$\rho\in T^*\cap \omega^{\textstyle n^\#+1}$, $\eta\hat{\
}l\subseteq\rho$. Clearly $\rho\in\bigcap_{\bar{\bf a}\in H_\rho} 
s_{\bar{\bf a}}$, so it is enough to show that for some $\rho$ the family
$H_\rho$ satisfies {\em (ii)}. Suppose that for each $\rho\in T^*\cap 
\omega^{\textstyle n^\#+1}$, $\rho\supseteq\eta\hat{\ }l$ we can find
$k_\rho<k^\#$ and $c_\rho$ such that
$c_\rho\in\bigcap\{\bar{\bf a}(k_\rho): \bar{\bf a}\in H_\rho\}$. Put
\[\bar{\bf a}^*(k)=f(\rho_k)\backslash\{c_\rho:\rho\in T^*\cap 
\omega^{\textstyle n^\#+1}\ \&\ \eta\hat{\ }l\subseteq\rho \}.\]
Let $\rho^+\in T^*\cap \omega^{\textstyle n^\#}$ be such that 
\[f(\rho^+)=\max\{f(\rho): \rho\in T^*\cap \omega^{\textstyle n^\#},
\eta\hat{\ }l \subseteq\rho\}.\]
By condition ($\gamma$) we get
\[\rest\{\rho\in T^*\cap \omega^{\textstyle n^\#+1}: \eta\hat{\
}l\subseteq\rho\}\rest \leq \prod\{f(\nu): f(\nu)\leq f(\rho^+)\}.\] 
Now, for each $k<k^\#$ we have $f(\rho^+)<f(\rho_k)$ (recall that
$\eta\hat{\ }l\subseteq \rho^+$, $\eta\hat{\ }m\subseteq\rho_k$ and
$l<m$ so condition ($\gamma$) works). Hence
\[\nor_{\rho_k}(\bar{\bf a}^*(k))\geq\frac{g(\rho_k)}{\prod\{f(\nu):
f(\nu)\leq f(\rho^+)\}} \geq \frac{g(\rho_k)}{\prod\{f(\nu):f(\nu)
<f(\rho_k)\}} > n^\#\]
Thus $\bar{\bf a}^*\in H$. Since $c_\rho\notin\bar{\bf a}^*(k_\rho)$ we have
$\bar{\bf a}^*\notin H_\rho$ for every $\rho\in T^*\cap \omega^{\textstyle
n^\#+1}$, $\eta\hat{\ }l\subseteq\rho$. Since $\bigcup\{H_\rho: \rho\in
T^*\cap \omega^{\textstyle n^\#+1}, \eta\hat{\ }l\subseteq\rho\}= H$ we
get a contradiction. The claim is proved.
\smallskip

Now let $H'\subseteq H$ be a family given by the claim. Condition {\em
(ii)} implies that 
\[r\forces_{\p}\mbox{``the family } \{T_{i_{\bar{\bf a}}}: \bar{\bf a}\in H'\}
\mbox{ has no upper bound in } \p_2 \mbox{''.}\]
Since $\rest H\rest\leq\prod\{2^{f(\rho_k)}:k<k^\#\}$ we have that for 
each $\rho\in T^*\cap\omega^{\textstyle n^\# + 1}$
\[\rest H'\rest\leq\prod\{2^{f(\nu)}:f(\nu)<f(\rho)\}.\]
Hence we may apply \ref{norm}~2) to conclude that for every 
$\rho\supseteq\eta\hat{\ }l$ , $\lh\rho\geq n^\#+1$:
\[r\forces_{\p}\mbox{``if }\rho\in\bigcap_{\bar{\bf a}\in H'}
\dot{S}_{i_{\bar{\bf a}}}\mbox{ then }\bigcap_{\bar{\bf a}\in H'}\suc_{
\dot{S}_{i_{\bar{\bf a}}}}(\rho)\neq\emptyset\mbox{''.}\]
Thus
\[r\forces_\p\mbox{``the family }\{\dot{S}_{i_{\bar{\bf a}}}:\bar{\bf
a}\in H'\} 
\mbox{ has an upper bound''.}\]
Since $r\forces_\p\mbox{``}\dot{S}_{i_{\bar{\bf a}}}\forces_{\p_2} 
\dot{T}_{i_{\bar{\bf a}}}\in\dot{R}_{n^*}\mbox{''}$ we get a
contradiction. \QED 

\noindent{\sc Remark:}\hspace{0.15in} 1) In the above theorem we worked in 
the model $\V[G]$ for technical reasons only. The assertion of the theorem 
can be proved in ZFC.

2) The forcing notion $\p_2$ is a special case of the forcing studied in 
[Sh1].

\begin{problem}
\hspace{0.15in}Does there exist a ccc Souslin forcing $\p$ such that

\noindent 1)\ \ \ $\p$ is homegeneous (i.e. for each $p\in\p$,
$\forces_{\p}$``there exists a generic filter $G$ over $\V$ such that
$p\in G$'')

\noindent 2)\ \ \ $\not\forces_{\p}$``$\hat{\p}$ is $\sigma$-centered''?
\end{problem}

\section{On ``small subsets of $\p$ are $\sigma$-centered''.}

Our next example is connected with the following, still open,
question:
\begin{problem}
\hspace{0.15in}Assume that for each ccc Souslin forcing $\p$ every
set $Q\!\in\! [\p]^{\textstyle \omega_{1}}$ is $\sigma$-centered (in
$\p$).

Does ${\bf MA}_{\omega_{1}}(Souslin)$ hold true?
\end{problem}
As an illustration of this subject let us recall a property of Random
(Solovay) Algebra $\B$ (see [BaJ]):
\begin{quotation}
\noindent if every $B\in [\B]^{\omega_{1}}$ is $\sigma$-centered 

\noindent then the real line can not be covered by $\omega_{1}$ null
sets and consequently ${\bf MA}_{\omega_{1}}(\B)$ holds true.
\end{quotation}
Our example shows that the above property of the algebra $\B$ does not
extend for other forcing notions. Let 
\[\p_3=\{(n,T): n\in \omega \ \& \ T\subseteq 2^{\textstyle
<\!\omega}\mbox{ is a tree } \ \& \ (\forall t\!\in\!
T\cap 2^{\textstyle n})(\mu([T_{t}])>0)\}.\] 
The order is defined by
\begin{quotation}
$(n_1,T_1)\leq(n_2,T_2)$ if and only if

$n_1\leq n_2$, $T_2\subseteq T_1$ and $T_1\rest n_1 = T_2\rest n_1$.
\end{quotation}

\begin{lemma}
\hspace{0.15in} $\p_3$ is a $\sigma$-linked Souslin forcing which is
not $\sigma$-centered.
\end{lemma}

\Proof Note that $\forces_{\p_3}$``there exists a perfect set of random
reals over $\V$''. Hence $\p_3$ is not $\sigma$-centered. To show that
it is $\sigma$-linked define sets $U(W,n,m)$ for $n<m<\omega$ and
finite trees $W\subseteq 2^{\textstyle \leq\! m}$:
\[U(W,n,m)=\{(n,T)\!\in\!\p_3: T\rest m\!=\!W\ \&\ (\forall t\in
T\!\cap\!2^n)(\mu([T_t]) > W(t)/2^{m+1})\}\]
where $W(t)=\rest\{s\in W\!\cap\! 2^m: t\subseteq s\}\rest$ (for $t\in
W\!\cap\!2^n$). Clearly each set $U(W,n,m)$ is linked (i.e. each two
members of it are compatible in $\p_3$) and $\p_3=\bigcup\{U(W,n,m):
n<m<\omega\ \&\ W\subseteq 2^{\textstyle \leq\! m}\}$. Since obviously
$\p_3$ is Souslin we are done. \QED

Let $\B(\kappa)$ stand for Random Algebra for adding $\kappa$ many
random reals. This is the measure algebra of the space $2^{\textstyle
\kappa}$. 

\begin{theorem}
\hspace{0.15in} Assume $\V\models$ {\bf CH}. Let
$G\subseteq\B(\omega_2)$ be a generic set over $\V$. Then, in $\V[G]$
\begin{description}
\item[(i) ] Martin axiom fails for $\p_3$ but
\item[(ii) ] each $Q\in [\p_3]^{\omega_{1}}$ is $\sigma$-centered (in $\p_3$).
\end{description}
\end{theorem}

\Proof Cichon proved that one random real does not produce a perfect
set of random reals (see [BaJ]). Hence in $\V[G]$ there is no perfect
set of random reals over $\V$. Consequently the first assertion is
satisfied in $\V[G]$. Since $\V[G]\models$ ``each
$B\in[\B]^{\textstyle \omega_{1}}$ is $\sigma$-centered in $\B$''
(compare section~\ref{notsigma}) it is enough to show the following
\smallskip

\noindent {\sc Claim:}\hspace{0.15in} Suppose that each
$B\in[\B]^{\textstyle \omega_{1}}$ is $\sigma$-centered. Then every
set $Q\in[\p_3]^{\textstyle \omega_{1}}$ is $\sigma$-centered.
\smallskip

\noindent Indeed, let $Q\in[\p_3]^{\textstyle \omega_{1}}$. For $n\in
\omega$ and $t\in 2^{\textstyle n}$ put
\[B(t,n)=\{[T_t]: (n,T)\in Q \ \&\ t\in T\}.\]
By our assumption we find sets $B(t,n,k)$ for $k,n\in \omega$, $t\in
2^{\textstyle n}$ such that $B(t,n)=\bigcup_{k\in \omega}B(t,n,k)$ and
for each $A_1 , A_2\in B(t,n,k)$ the set $A_1 \cap A_2$ is of positive
measure. Now define sets $Q(n,W,\sigma)$ for $n\in \omega$, a finite
tree $W\subseteq 2^{\textstyle \leq n}$ and a function $\sigma: W\cap
2^{\textstyle n} \longrightarrow  \omega$:
\[Q(n,W,\sigma)=\{(n,T)\in Q: T\rest n=W\ \&\ (\forall t\in T\cap
2^{\textstyle n})([T_t]\in B(t,n,\sigma(t)))\}.\] 
Note that if $(n,T_1),(n,T_2)\in Q(n,W,\sigma)$ then for each $t\in
W\cap 2^{\textstyle n}$ the set $[(T_1)_t]\cap[(T_2)_t]$ is of
positive measure. Consequently each $Q(n,W,\sigma)$ is linked and we
are done. \QED

\section{A $\sigma$-centered example}

In this section we define a very simple $\sigma$-centered Souslin forcing 
notion. Next we show that in any generic extension of some model of {\bf CH} 
via finite support iteration of the Dominating (Hechler) Algebra, Martin Axiom
fails for this forcing notion. Consequently we get the consistency of
the following sentence:
\begin{quotation}
{\em \noindent any union of less than continuum meager sets is meager
+ $\neg${\bf CH} + {\bf MA} fails for some $\sigma$-centered Souslin
forcing.} 
\end{quotation}

Our example $\p_4$ consists of all pairs $(n,F)$ such that $n\in
\omega$, $F\in [\can]^{\textstyle <\omega}$ and all elements of the
list $\{x\rest n: x\in F\}$ are distinct. $\p_4$ is ordered by
\begin{quotation}
\noindent $(n,F)\leq(n',F')$ if and only if

\noindent $n\leq n'$, $F\subseteq F'$ and $\{x\rest n: x\in
F\}=\{x\rest n: x\in F'\}$.
\end{quotation}

\begin{lemma}
\hspace{0.15in}$\p_4$ is a $\sigma$-centered Souslin forcing.
\end{lemma}

\Proof Clearly $\p_4$ is Souslin (even Borel). To show that $\p_4$ is
$\sigma$-centered note that if $\{x\rest n: x\in F_0\}=\ldots=\{x\rest
n: x\in F_k\}$ then the conditions $(n,F_0 ),\ldots,(n,F_k )$ are
compatible (if $m$ is large enough then $(m,F_0\cup\ldots\cup F_k )$
is a witness for this). \QED

Now we want to define the model we will start with. At the beginnig we
work in $\L$. Applying the technology of [Sh] we can construct a
sequence $(\p_{\xi}:\xi\leq \omega_{1})$ of forcing notions such that
for each $\alpha,\beta<\omega_{1}$, $\xi\leq \omega_{1}$:
\begin{description}
\item[(1) ] if $\alpha<\beta$ then $\p_{\alpha}$ is a complete suborder
of $\p_{\beta}$,
\item[(2) ] there is $\gamma>\beta$ such that
$\p_{\gamma+1}=\p_{\gamma}*\dot{\D}_{\alpha}$, where
$\dot{\D}_{\alpha}$ is the $\p_{\gamma}$-name for finite support,
$\alpha$ in length, iteration of Hechler forcing,
\item[(3) ] $\p_{\xi}$ satisfies ccc,
\item[(4) ] if $\xi$ is limit then
$\p_{\xi}=\stackrel{\longrightarrow}{\lim}_{\zeta<\xi}\p_{\zeta}$,
\item[(5) ] $\p_{\omega_{1}}\forces$ ``every projective set of reals
has Baire property''
\end{description}
(for details see also [JR]). Recall that Hechler forcing $\D$ consits
of all pairs $(n,f)$ such that $n\in \omega$, $f\in\baire$. These pairs
are ordered by
\begin{quotation}
\noindent $(n,f)\leq(n',f')$ if and only if

\noindent $n\leq n'$, $f\rest n=f'\rest n'$ and $f(k)\leq f'(k)$ for
all $k\in \omega$.
\end{quotation}

Suppose $G\subseteq\p_{\omega_{1}}$ is a generic set over $\L$. We
work in $\L[G]$. For distinct $x,y\in\can$ we define $h(x,y)=\min\{n:
x(n)\neq y(n)\}$. Easy calculations show the following
\begin{lemma}
\label{Raisonnier}
\hspace{0.15in} Let $b\subseteq \omega$. Then the following conditions
are equivalent:
\smallskip

\noindent {\bf (i)} there exists a {\em Borel} equivalence relation $R$ on
$\can$ with countable many equivalence classes such that 
$\{h(x,y): x,y\in\can\cap\L \ \&\  x\neq y \ \&\  R(x,y)\}\subseteq b$,
\smallskip

\noindent {\bf (ii)} there exists an equivalence relation $R$ on
$\can$ with countable many equivalence classes such that 
$\{h(x,y): x,y\in\can\cap\L \ \&\  x\!\neq\! y \ \&\  R(x,y)\}\subseteq b$,
\smallskip

\noindent {\bf (iii)} there exist sets $Y_n\subseteq\can$ (for $n\in \omega$)
such that $\L\cap\can\subseteq\bigcup_{n\in \omega}Y_n$ and
$\bigcup_{n\in \omega}\{h(x,y):x\!\neq\! y \  \&\  x,y\in Y_n\}\subseteq b$,
\smallskip

\noindent {\bf (iv)} $(\exists f\!:\!2^{\textstyle <\omega}\!\rightarrow\!
2)(\forall x\!\in\!\can\!\cap\!\L)(\exists m\!\in\!\omega)(\forall
n\!>\!m)(n\notin b \Rightarrow f(x\rest n)=x(n))$.\hspace{-0.2in}\QED 

\end{lemma}

The Raisonnier filter $\cal F$ consists of all sets $b\subseteq\omega$
satisfying one of the conditions of~\ref{Raisonnier} (cf [Ra]). $\cal F$
is a proper filter on $\omega$. Directly from {\bf (iv)}
of~\ref{Raisonnier} one can see that $\cal F$ is a $\Sigma^1_3$-subset
of $\can$. Consequently it has Baire property (recall that we are in
$\L[G]$). 

\begin{theorem}
{\em (Talagrand, [Ta])}\hspace{0.15in} For any proper filter $F$ on
$\omega$ the following conditions are equivalent:
\vspace{0.1in}

\noindent {\bf (i)} $F$ does not have Baire property,
\vspace{0.1in}

\noindent {\bf (ii)} for every increasing sequence $(n_k :k\in
\omega)$ of integers there exists $b\in F$ such that
$(\exists^{\infty}_k)(b\cap[n_k,n_{k+1})=\emptyset)$. \QED
\end{theorem}
Applying the above theorem we can find an increasing function
$r\in\baire\cap\L[G]$ such that (in $\L[G]$)
$$(\forall b\in{\cal F})(\forall^{\infty}_{k})(b\cap [r(k),r(k+1))
\neq\emptyset)$$
Let $\dot{r}$ be the $\p_{\omega_{1}}$-name for $r$ and let
$\alpha_0<\omega_{1}$ be such that $\dot{r}$ is a $\p_{\alpha_0}$-name. 

Our basic model will be $\L[r]$

\begin{theorem}
\hspace{0.15in} Let $\kappa$ be a regular cardinal. Let $\D_\kappa$ be
the finite support iteration of Hechler forcing of the length
$\kappa$. Suppose $H\subseteq\D_\kappa$ is a generic set over $\L[r]$.
Then
\begin{quotation}
$\L[r][H]\models$``there is no $\p_4$-generic over $\L[r]$.''
\end{quotation}
\end{theorem}

\Proof Assume not. Let $H^*\in\L[r][H]$ be a $\p_4$-generic over
$\L[r]$. Put $T=\bigcup\{F: (\exists n\!\in\!\omega)((n,F)\in H^*)\}$.
Then in $\L[r][H^*]$ we have:
\begin{description}
\item[(6) ] $T$ is a closed subset of $\can$,

\item[(7) ] $(\exists^\infty_k)(\forall x,y\in T)(x\neq y\Rightarrow
h(x,y)\notin [r(k),r(k+1)))$ and

\item[(8) ]  $(\forall x\in\can\cap\L)(\exists q\in Q)(q+x\in T)$
\end{description}
\noindent (Q stands for the set of all sequences eventually equal 0, +
denotes the addition modulo 2). Since both (7) and (8) are absolute
($\Pi^1_2$) sentences they are satisfied in $\L[r][H]$ too. Let
$\dot{T}\in\L[r]$ be a $\p_{\kappa}$-name for $T$. Since $T$ is a
closed subset of $\can$ we can think of $\dot{T}$ as a name for a
real.

Now we work in $\L[r]$. Let $p\in\D_{\kappa}\cap \L[r]$ be such that
\begin{quotation}
$p\forces$``$\dot{T}$ satisfies (6), (7) and (8)''.
\end{quotation}
By Souslin forcing properties (see \S 1 of [JS1]) we find a (closed)
countable set $S\subseteq\kappa$ such that:
\begin{description}
\item[(9) ] $\dot{T}$ is a $\D_{\kappa}\rest S$-name,
$p\in\D_{\kappa}\rest S$ and

\item[(10) ] $\D_{\kappa}\rest S$ is a complete suborder of $\D_{\kappa}$.
\end{description}
\noindent Since (6)-(8) are absolute we get
\begin{description}
\item[(11) ] $p\forces_{\D_{\kappa}\rest S}$``$\dot{T}$ satisfies (6),
(7) and (8)''.
\end{description}
But $\D_{\kappa}\rest S$ is isomorphic to finite support iteration of
Hechler forcing of the countable length $\alpha$
($\alpha<\omega_{1}$). Thus we can treat $\dot{T}$ as a
$\D_{\alpha}$-name and $p$ as a condition in $\D_{\alpha}$. Then, in
$\L[r]$
\begin{description}
\item[(12) ] $p\forces_{\D_{\alpha}\rest S}$``$\dot{T}$ satisfies (6),
(7) and (8)''.
\end{description}
By (2) we find $\beta>\alpha_0$ such that
\begin{description}
\item[(13) ] $\p_{\gamma+1}=\p_{\gamma}*\dot{\D_{\alpha}}$ and
\item[(14) ] $p\equiv({\bf 1},p)$ interpreted as a member of
$\p_{\gamma+1}$ belongs to $G$.
\end{description}
By Souslin forcing properties (12) holds true in $\L[G\cap\p_{\gamma}]$
and hence
\begin{description}
\item[(15) ] $\L[G\cap\p_{\gamma+1}]\models$``$\dot{T}^G$ satisfies
(6), (7) and (8)''
\end{description}
(we treat here $\dot{T}$ as a $\p_{\gamma+1}$-name). Let $b=\{h(x,y):
x\neq y\ \& \ x,y\in\dot{T}^G\}\in[\omega]^{\textstyle <\omega}\cap\L[G]$. By
(15) and by Shoenfield absoluteness we have
\begin{description}
\item[(16) ] $\L[G]\models$``$\dot{T}^G$ satisfies (6), (7) and (8)''.
\end{description}
Since $\{h(x,y): x\neq y\ \& \ x,y\in\dot{T}^G\}=\{h(x,y): x\neq y\ \&
\ x,y\in\dot{T}^G+q \}$ we conclude that
\begin{description}
\item[(17) ] $\L[G]\models$``sets $\dot{T}^G + q$ (for $q\in Q$)
witness that $b\in{\cal F}$'' and
\item[(18) ] $\L[G]\models\
(\exists^\infty_k)(b\cap[r(k),r(k+1))=\emptyset)$. 
\end{description}
The last condition contradicts our choice of $r$. \QED

Since $\forces_{\D_{\kappa}}$``any union of less than $\kappa$ meager
sets is meager'' we get
\begin{corollary}
\hspace{0.1in} The following theory is consistent:

\noindent {\bf ZFC} + $\neg${\bf CH} + ``Martin Axiom fails for some
$\sigma$-centered Souslin forcing'' + ``any union of less than
continuum meager sets is meager''. \QED
\end{corollary}

\section{On Souslin not ccc}

In this section we will give a negative answer to the following question of
Woodin:
\begin{quotation}
\noindent If $\p$ is a Souslin forcing notion which is not ccc

\noindent then there exists a perfect set $T\subseteq\p$ such that
each distinct $t_{1}, 
t_{2}\in T$ are incompatible.
\end{quotation}

\noindent Recall that in the case of non-ccc partial orders we do not
require Souslin forcings to satisfy the condition:

 ``the set $\{(p,q): p \mbox{ is incompatible with } q\}$ is
$\Sigma^1_1$''.  

\noindent Thus a forcing notion $\p$ is {\em is Souslin not ccc} if
both $\p$ and $\leq_{\p}$ are analytic sets. The reason for this is
that we want to cover in our definition various standard forcing
notions with simple definitions for which incompatibility is not
analytic (e.g. Laver forcing).
\medskip
 
Let $\Q$ be the following partially ordered set:
\medskip

$W\in\Q$ if $W$ is a finite set of pairs $(\alpha,\beta)$, $\alpha\leq\beta
<\omega_{1}$ such that 
if $(\alpha_{1},\beta_{1}),(\alpha_{2},\beta_{2})$ are in $W$, then $\beta_{1}
<\alpha_{2}$ or $\beta_{2}<\alpha_{1}$.
 
$\Q$ is ordered by the inclusion.
\medskip

It follows from [Je1] that $\Q$ is proper. Clearly
$\rest\Q\rest=\omega_{1}$. 
Next define a forcing notion $\p_5$. It consists of all $r\in\baire$
such that 
$r$ codes a pair $(E^{r},w^{r})$ where
\begin{enumerate}
\item $E^{r}$ is a relation on $\omega$ such that 
$(\omega,E^r)\models{\bf ZFC}^{-}$ and $E^r$ encodes all elements of
$\omega\cup\{\omega\}$. 
\item $w^{r}\in \omega$ and $E^{r}\models$``$w^{r}\in\Q$''. 
\end{enumerate}

We say that a one-to-one function $f\in\baire$ {\em interprets}
$E^{r_1}$ {\em in} $E^{r_2}$ if there exists $n\in \omega$ such that
rng$(f)=\{k\in \omega:E^{r_2}(k,n)\}$ and $E^{r_1}(l,k)\equiv
E^{r_2}(f(l),f(k))$. 

If $f$ interprets $E^{r_1}$ in $E^{r_2}$ then $E^{r_2}$ may ``discover'' 
that some of the ordinals of $E^{r_{1}}$ are not
ordinals (i.e. not well-founded). Let
$w(r_{1},r_{2},f)=w^{r_{1}}\cap\{(\alpha, \beta):
\alpha\leq\beta\mbox{ are ordinals in } E^{r_{2}}\}$. Then, in $E^{r_2}$,  
$w(r_{1},r_{2},f)$ is an initial segment of $w^{r_{1}}$ and it is in 
$\Q^{E^{r_{2}}}$.

Now we can define the order $\leq$ on $\p_5$:
\begin{quotation}
\noindent $r_{1}\leq r_{2}$ if and only if $r_{1}=r_{2}$ or

\noindent {\bf there exists} $f\in\baire$ which interprets $E^{r_1}$
in $E^{r_2}$ and such that
\[(\omega,E^{r_2})\models w(r_1,r_2,f)\subseteq w^{r_2}.\] 
\end{quotation}

Obviously both $\p_5$ and the order $\leq$ are $\Sigma^1_1$-sets. 

For $r\in\p_5$ we define $W(r)$ as
$w^r\cap\{(\alpha,\beta):\alpha\leq\beta\mbox{ are well founded }\}$.
Note that $W(r_1)=W(r_2)$ implies $r_1$ and $r_2$ are equivalent in
$\p_5$ (i.e they have the same compatible elements of $\p_5$).
Consequently $\Q$ may be densely embedded into the complete Boolean
algebra determined by $\p_5$. 
It follows from [Je1] that $\p$ is 
proper, it is Souslin and it does not satisfy the countable chain condition.
Moreover, if $\omega_{1}<\can$ then $\p$ does not contain a perfect set
of pairwise incompatible elements (recall $\rest\Q\rest=\omega_{1}$). \QED

An interesting question appears here:
\begin{quotation}
\noindent Suppose $\p$ is $\omega$-proper and Souslin. 

\noindent Does there exists a perfect set of pairwise incompatible
elements of $\p$? 
\end{quotation}
The negative answer to this question is given by the following result.
\begin{theorem}
\hspace{0.15in}Assume $\omega_1<\mbox{cf}(\can)$. There exists an $\omega$-proper
Souslin not ccc forcing notion $\p_5^*$ with no perfect set of
pairwise incompatible elements. 
\end{theorem}

\Proof Let $\delta\leq\omega_1$ be additively indecomposable. Let
$\Q^*$ be the order defined by:
\begin{quotation}
\noindent $W\in\Q^*$ if and only if 

\noindent $W$ is a countable set of pairs $(\alpha,\beta)$,
$\alpha\leq\beta<\omega_1$ such that
\begin{itemize}
\item $(\alpha_1,\beta_1),(\alpha_2,\beta_2)\in W\Rightarrow
\beta_1<\alpha_2\mbox{ or }\beta_2<\alpha_1$,

\item $\{(\alpha,\beta)\in W:\alpha\neq\beta\}$ is finite,

\item the order type of the set $\{\alpha:(\exists\beta)((\alpha,\beta)\in
W)\}$ is less than $\delta$.
\end{itemize}
\noindent $\Q^*$ is ordered by the inclusion.
\end{quotation}
 
It follows from Chapter XVII, \S 3 of [Sh 2] that $\Q^*$ is
$\alpha$-proper for each $\alpha<\omega_1$.

Now we can repeat the coding procedure that we applied to define the
forcing notion $\p_5$. Thus we get the Souslin forcing notion $\p_5^*$
such that $\Q^*$ can be densely embedded in the Boolean algebra
determined by $\p_5^*$. 

For $W\in\Q^*$ let heart$(W)=\{(\alpha,\beta)\in W:\alpha\neq\beta\}$.

Assume that $\{(E^{r_\eta},w^{r_\eta}:\eta\in\can\}\subseteq\p_5^*$
is a perfect set of pairwise incompatible elements. Let $W_\eta$ be
the well-founded part of $w^{r_\eta}$. Since w.l.o.g we can assume
that $\sup\{\beta:(\exists \alpha)((\alpha,\beta)\in W_\eta)\}$ is
constant and heart$(W_\eta)$ is constant we easily get a
contradiction.\QED

\section{On ccc $\Sigma^1_2$}

Souslin ccc notions of forcing are indestructible ccc (see [JS1]):
\begin{quotation}
\noindent {\em Suppose $\p$ is a ccc Souslin notion of forcing. Let
$\Q$ be a ccc 
forcing notion. Then $\forces_{\Q}$``\ $\hat{\p}$ is ccc\ ''.}
\end{quotation}
The above property does not hold true for more complicated forcing
notions. In this section we show that there may exist two ccc
$\Sigma^1_2$-notions of forcing $\p_6$ and $\p_6^{*}$ such that
$\p_6\times\p_6^{*}$ does not satisfy ccc.

We start with $\V=\L$. Let $\Q$ be a ccc notion of forcing such that
\[\forces_{\Q}\ {\bf MA} + \neg{\bf CH}\]
Let $G\subseteq\Q$ be a generic set over $\L$ and let $r$ be a random
real over $\L[G]$. Recall that by theorem of Roitman (cf [Ro]) we have
$\L[G][r]\models {\bf MA}(\sigma\mbox{-centered})$. 

Fix a sequence $(f_\alpha:\alpha<\omega_{1})\in\L$ of one-to-one
functions
$f_\alpha:\alpha\stackrel{\mbox{1--1}}{\longrightarrow}\omega$ and
define in $\L[r]$ sets $E_1, E_2$ by
\[E_i=\{\{\alpha,\beta\}\in[\omega_{1}]^2: \beta<\alpha\ \&\
r(f_\alpha(\beta))=i\} \ \ \ \mbox{for } i=0,1.\]
We define forcing notions $\p_6, \p_6^*$:
\[\p_6=\{H\in[\omega_{1}]^2: [H]^2\subseteq E_0\},\]
\[\p_6^*=\{H\in[\omega_{1}]^2: [H]^2\subseteq E_1\}.\]
Orders are inclusions. 

\noindent Both $\p_6$ and $\p_6^*$ are elements of $\L[r]$. Moreover
they can be thought of as subsets of $\L[r]\cap\can$. Applying {\bf
MA}($\sigma$-centered) we get that (cf [Je]):
\[\L[G][r]\models\mbox{``any subset of } \L[r]\cap\can \mbox{ is a
relative }\Sigma^0_2\mbox{-set''}.\]
Consequently
\[\L[G][r]\models\mbox{``any subset of } \L[r]\cap\can \mbox{ is
$\Sigma^1_2$''}.\] 
Thus $\p_6$ and $\p_6^*$ are $\Sigma^1_2$-notions of forcing in
$\L[G][r]$ (i.e. both $\p_6$, $\p_6^*$ and orders and the relations of
incompatibility are $\Sigma^1_2$-sets). Roitman proved the following
\begin{theorem}
{\em (Roitman, Prop.4.6 of [Ro])}\hspace{0.15in}In $\L[G][r]$ both $\p_6$
and $\p_6^*$ satisfy ccc and $\p_6\times\p_6^*$ does not satisfy ccc. \QED
\end{theorem}
\begin{corollary}
\hspace{0.15in}The following theory is consistent:

\noindent {\bf ZFC} + {\bf MA}($\sigma$-centered) + $\neg${\bf CH} +
``there exist ccc $\Sigma^1_2$-notions of forcing $\p_6$, $\p_6^*$
such that $\p_6^*\forces$`` $\hat{\p}_6$ is not ccc''.\QED
\end{corollary}

\begin{problem}
\hspace{0.15in} Is there a ccc Souslin forcing notion $\p$ such that
MA($\p$) always fails after adding a random real?
\end{problem}

\eject

\bigskip
\bigskip
\bigskip
\bigskip
\bigskip

{\sc References}

[Ba] J.Baumgartner, {\em Iterated forcing} in Surveys in Set Theory, ed. by
A.R.D.Mathias, London Math. Soc., Lecture Notes 87.

[BJ] J.Bagaria, H.Judah, {\em Amoeba forcing, Souslin absoluteness and 
additivity of measure}, Proceedings of the 1989 MSRI Workshop on Set
Theory of the Reals (to appear).

[BaJ] T.Bartoszynski, H.Judah, {\em Jumping with random reals},

[Je] T.Jech, {\em Set Theory}, Academic Press, New York 1978.

[Je1] T.Jech, {\em Multiple Forcing}, Cambridge Tracts in Mathematics
88, Cambridge 1986.

[JS1] H.Judah, S.Shelah, {\em Souslin forcing}, Journal of Symbolic Logic, 
53(1988).

[JS2] H.Judah, S.Shelah, {\em Martin's axioms, measurability and 
equiconsistency results}, Journal of Symbolic Logic, 54(1989).

[JR] H.Judah, A.Roslanowski, {\em On Shelah's amalgamation},
Proceedings of the Winter Institute on Set Theory of the Reals,
Bar-Ilan 1991 (to appear).

[Ra] J.Raisonnier, {\em A mathematical proof of Shelah's theorem},
Israel Journal of Mathematics 48(1984).

[Ro] J.Roitman, {\em Adding a random or a Cohen real: topological
consequences and the effect of Martin's axiom}, Fundamenta
Mathematicae vol CIII (1979).

[Sh] S.Shelah, {\em Can you take Solovay's inaccessible away?}, Israel
Journal of Mathematics, 48(1984).

[Sh1] S.Shelah, ???{\em Vive la Difference I (new version of models with no 
isomorphic ultra powers}, Proceedings of the Conference in Set Theory, 
MSRI 10/89,

[Sh2] S.Shelah, {\em Proper and Improper Forcing}, in preparation.

[Ta] M.Talagrand, {\em Compacts de functions mesurables et filtres non
mesurables}, Studia Mathematica, 67(1980).

[To] S.Todorcevic, {\em Two examples of Borel partially ordered sets
with the countable chain conditions}, Proceedings of the American
Mathematical Society, 112(1991).

\end{document}